%% file: Zhukov_eng18Arxiv.tex
\documentclass[11pt]{article}

\usepackage[utf8]{inputenc}
\usepackage[russian,english]{babel}
\usepackage{amssymb,amsfonts,amsmath,mathtext,mathrsfs}
\usepackage{tikz}

\def\cB{{\cal B}}

\def\cB{{\cal B}}

\def\cF{{\cal F}}

\def\cG{{\mathscr{G}}}

\def\cL{{\mathscr{L}}}

\def\F{{\mathbb F}}

\usepackage{color}

\definecolor{MyGrayc}{gray}{0.30}
\definecolor{MyGreen}{rgb}{0.0,0.39368872549019607,0.0}
\definecolor{MyRed}{rgb}{0.8, 0, 0}
\definecolor{MyBlue}{rgb}{0,0,0.8}
\definecolor{my-color-chocolate}{rgb}{0.8,0.5,0.4}
\definecolor{my-color-dark-magenta}{rgb}{0.5472273284313726,0.0,0.5472273284313726}

\newcommand{\warn}{% \color{MyBlue}
}

\newcounter{todomCounter}

\newenvironment{Proof}{\par\noindent{\bf Proof.}}
{\hfill$\scriptstyle\blacksquare$}

\usepackage{theorem}

\newtheorem{Th}{Theorem}[section]
\newtheorem{Prop}{proposition}[section]
\newtheorem{Le}[Prop]{Lemma}

{\theorembodyfont{\rmfamily}
  \newtheorem{Def}[Prop]{Definition}
  \newtheorem{Ex}[Prop]{Example}
}
\newtheorem{corollary}[Prop]{Corollary}

\usepackage{soul}

\usepackage[normalem]{ulem}

\newcommand{\comment}[1]{}
\title{Lagrangian subspaces, delta-matroids and four-term relations\footnote{The article was prepared within the framework of the Academic Fund Program at the National Research University Higher School of Economics (HSE)
in 2016|2017 (grant 16-05-0007) and supported within the framework of a
subsidy granted to the HSE by the Government of the Russian Federation
for the implementation of the Global Competitiveness Program.}}

\author{V.I.~Zhukov% \thanks{National Research University Higher School of Economics; e-mail: slava.zhukov@list.ru}
}
\date{}

\begin{document}
\maketitle

\begin{abstract}
  Finite order invariants (Vassiliev invariants) of knots are expressed
  in terms of weight systems, that is, functions on chord diagrams satisfying the four-term relations.
  Weight systems have graph analogues, so-called $4$-invariants of graphs,
  i.e. functions on graphs that satisfy the four-term relations for graphs.
  Each $4$-invariant determines a weight system.

  The notion of weight
  system \comment{\st{naturally generalizes}} is naturally generalized for the case of embedded graphs with an
  arbitrary number of vertices.  Such embedded graphs correspond to links;
  to each component of a link there corresponds a vertex of an embedded graph.
  Recently, two approaches have been suggested
  to extend the notion of $4$-invariants of  graphs to the case
  of combinatorial structures corresponding to  embedded graphs with an arbitrary number of
  vertices.  The first approach is due to V.~Kleptsyn and E.~Smirnov, who considered
  functions on Lagrangian subspaces in a $2n$-dimensional space over $\mathbb{F}_2$ endowed with a standard
  symplectic form % {\warn \st{Symplectic}} \st{spaces over the field of two elements}
  and introduced
  four-term relations for them.  On the other hand, the second approach, the one due to
  Zhukov and Lando,
  suggests four-term relations for functions on binary delta-matroids.  In this paper, we prove
  that the two approaches are equivalent.
\end{abstract}
%% ====

Finite order invariants
(Vassiliev invariants) of knots are expressed in terms of weight
systems, that is, functions on chord diagrams satisfying four-term relations.  The vector space over $\mathbb{C}$  spanned by chord diagrams
considered modulo
four-term relations is
supplied with a Hopf algebra structure.
The notion of weight system is naturally
extended from functions on chord diagrams
(which can be interpreted as embedded graphs with a single vertex) to
functions on arbitrary embedded graphs.

In~\cite{BGBS}, to each embedded graph a
Lagrangian subspace in a symplectic space over the field $\mathbb{F}_2$ is associated.  V.~Kleptsyn, E.~Smirnov
in~\cite{KS} rediscovered this construction. They introduced four-term relations
in the vector space spanned by Lagrangian
  subspaces,  and showed that linear functionals satisfying these
four-term relations produce weight systems.  They constructed a Hopf algebra of Lagrangian
subspaces and a quotient Hopf algebra of Lagrangian subspaces modulo the four-term relations.

Meanwhile, Lando and Zhukov in~\cite{LZ} constructed a Hopf algebra of binary delta-matroids,
introduced four-term relations for them and constructed a quotient Hopf algebra modulo the four-term relations.
The correspondence between delta-matroids and embedded graphs allows one to associate a weight system to a linear
functional on the latter Hopf algebra.  The main result of the present
paper is the proof of equivalence of these two approaches; in particular, we establish an
isomorphism between the Hopf algebra of Lagrangian subspaces and the Hopf algebra of binary
delta-matroids.  This isomorphism is given by the mapping $\nu_E$, which establishes (according to
  Theorem 2.1) a one-to-one correspondence between the set of Lagrangian subspaces in
$V_E$, the vector space spanned by the elements of a finite set $E$ as well as their duals,
and binary delta-matroids on the set $E$.

\section{Necessary information about delta-matroids}

A {\em set system\/} $(E;\Phi)$ is a pair consisting of a finite set~$E$ and a set $\Phi\subset 2^E$
of subsets of~$E$.
The set $E$ is called the {\em ground set} and the elements of the set $\Phi$ are called the {\em
feasible subsets}
of this system.

Two set systems
$(E_1; \Phi)$, $(E_2; \Phi_2)$ are said to be {\em isomorphic} if there exists a one-to-one correspondence
$E_1\rightarrow E_2$, which identifies the subsets $\Phi_1\subset 2^{E_1}$ with
$\Phi_2\subset 2^{E_2}$.  Below, we will not distinguish between isomorphic set systems.

A set system
$(E; \Phi)$ is said to be {\em proper} if the set $\Phi$ is nonempty.  In our paper, we consider only
proper set systems if otherwise is not stated explicitly.  We denote by $\Delta$ the set symmetric
difference operation, that is, $A\Delta B=(A\setminus B)\sqcup (B\setminus A)$.  A {\em delta-matroid} is a set system $(E; \Phi)$
that satisfies the following symmetric exchange axiom (SEA): {\em for any two feasible subsets $\phi_1$
  and $\phi_2\in\Phi$ and for any element $e\in\phi_1\Delta\phi_2$ there exists an element
  $e'\in\phi_1\Delta\phi_2$ such that $\phi_1\Delta\{e,e'\}\in\Phi$.}

Let $G$ be an (abstract) simple
graph. We will consider more general objects, namely, framed graphs,
that is, graphs all whose vertices are endowed with an element $0$ or $1$ of the field $\mathbb{F}_2$.  To each framed
graph $G$, with the set of vertices $V(G)$, one can assign its adjacency matrix $A(G)$
(of dimension $|V(G)|\times|V(G)|$) on the intersection of the row $v$ and the column $v'$ of which
($v\not=v'$), there is the element $1$ of the field $\mathbb{F}_2$ if the vertices $v$ and $v'$ are
neighbors (that is, are connected by an edge), and the element $0$, otherwise. In turn, the diagonal
elements are equal to the {frames} of the corresponding vertices.

A framed graph $G$ is said to
be {\em non-degenerate} if its adjacency matrix $A(G)$, considered as a matrix over the field
$\mathbb{F}_2$, is non-degenerate, i.e. if its determinant equals $1$.  Let us define the set system
$(V(G); \Phi(G))$, $\Phi(G)\subset 2^{V(G)}$ in the following way:
\begin{eqnarray*}
  V(G) &&\text{ is the set of the vertices of }G \\
  \Phi(G)&=&\{U\subset V(G)\ |\ G_{U} \text{ is non-degenerate} \},
\end{eqnarray*}
% where $V(G)$ is the set of the graph $G$ vertices and
% $\Phi(G)=\{U\subset V(G)| G_{U} \text{ is non-degenerate}\}$
where $G_{U}$ denotes the subgraph in~$G$  induced by the vertex set $U$.
\begin{Th}[\cite{Bo871}]\label{tDG}
  The set system $(V(G); \Phi(G))$ is a delta-matroid.
\end{Th}
We call this delta-matroid the {\em non-degeneracy delta-matroid} of the graph $G$.

Non-degeneracy delta-matroids of framed graphs are examples of binary delta-matroids.
To introduce the notion of binary delta-matroid, we need the operation of twisting.  For a set system
$D=(E;\Phi)$ and a subset $E'\subset E$, let us define the twist $D*E'$ of the set system $D$ by
the subset $E'$ by the equation
$$D*E'=(E;\Phi\Delta E') = (E;\{\phi\Delta E'|\phi\in\Phi\}).$$
Obviously, twisting of set systems by a subset is an involution, \mbox{$D*E'*E'=D$.}
\begin{Th}[\cite{Bo89}]
  The twist of a non-degeneracy delta-matroid of a framed graph by any subset is a
  delta-matroid.
\end{Th}
\begin{Def}[\cite{Bo89}]
  A {\em binary delta-matroid } is the result of twisting the non-degeneracy
    delta-matroid of a framed graph by (maybe an empty) subset.
\end{Def}
Denote by $\mathscr{B}_E$ the set of binary delta-matroids with the ground set $E$.

%% ================================================================
\section{Binary delta-matroids and Lagrangian subspaces (set-theoretic bijection)}
\label{sec:binary-delta-matr}

% \rus{В этом параграфе мы устанавливаем биективное соответствие между множеством бинарных дельта-матроидов на данном конечном множестве $E$ и множеством лагранжевых подпространств в симплектическом пространстве $V_E$ над
% полем $\mathbb{F}_2$ из двух элементов, ассоциированном с множеством $E$}
In this section we % give/
establish a one-to-one correspondence between the set of binary delta-matroids
(on a finite set $E$) and the set of Lagrangian subspaces in the symplectic space $V_E$ over the field
$\mathbb{F}_2$ associated with the set $E$.

Let $E$ be a finite set and $E^\vee$ be its copy.
Denote by $e^\vee$ the element of $E^\vee$ corresponding to the element $e$ in $E$.
We denote by $^\vee:E\sqcup E^\vee\rightarrow E\sqcup E^\vee$  the bijection of $E\sqcup E^\vee$,
which exchanges the elements $e$ and $e^\vee$ for all $e\in E$.
For $Y\subset E\sqcup E^\vee$, denote by $Y^\vee$ the image of $Y$ under the map $^\vee$.

A {\em symplectic structure} on a vector space is a nondegenerate skew symmetric form on it.
Symplectic structures exist only on even-dimensional spaces.
Denote by $V_E$ the $2|E|$-dimensional space over the field $\mathbb{F}_2$ spanned by the elements of the set $E\sqcup E^\vee$.
Let us introduce a symplectic structure  $(\cdot, \cdot)$
on $V_E$ %-- over $\mathbb{F}_2$
by the rule
$(e, e^\vee)=(e^\vee, e)=1$, and $(u, v)=0$ otherwise.

A subspace $L$ of a symplectic space is said to be {\em isotropic} if the restriction
 of the symplectic form to~$L$ is zero, i.e. $(u,v)=0$ for all $u$ and $v$
in $L$.
The dimension of an isotropic subspace  of a symplectic space cannot
exceed half of the dimension of the symplectic space itself.
An isotropic subspace whose dimension is half the dimension of the symplectic space
is called a {\em Lagrangian subspace}.
Denote by $\mathscr{L}_E$ the set of Lagrangian subspaces in $V_E$.

\begin{Def}(mapping $\nu_E$)\footnote{A similar mapping is considered in \cite{MG}.}
  \label{sec:def1} Let $L$ be an arbitrary Lagrangian subspace in $V_E$.
  Denote by $\nu_E(L)$ the set system $\nu_E(L)=(E; \Psi_L)$, where a subset $Y\subset E$
  belongs to $\Psi_L$ if and only if $L{\cap} \langle Y^\vee{\sqcup} (E{\setminus} Y)\rangle=0$;
  Here the angle brackets denote the vector subspace in $V_E$ spanned by the elements inside, and $0$ is the zero vector of the space $V_E$.

\end{Def}

\begin{Ex}
  Let $E$ be a $2$-element set, $E=\{1, 2\}$, then $L=\langle 1^\vee + 2 + 2^\vee, 1 + 2 \rangle$ is a Lagrangian subspace
 in $V_E$. It consists of four elements, namely, $0$, $1^\vee+2+2^\vee$, $1+2$, $1+1^\vee+2^\vee$. Then
 $\nu_E(L)=(E;\{\{1\}, \{2\}, \{1, 2\}\})$. (In [10], this set system is denoted by $s_{25}$.)
 Indeed, we have

 \begin{tabular}{lll}
   for $Y=\emptyset$, & $\langle Y^\vee \sqcup (E\setminus Y)\rangle=\langle{1, 2}\rangle$, & $L\cap \langle{1, 2}\rangle\ni 1+2 $; \\
   for $Y=\{1\}$, & $\langle Y^\vee \sqcup (E\setminus Y)\rangle=\langle{1^\vee, 2}\rangle$, & $L\cap \langle{1^\vee, 2}\rangle=0 $; \\
   for $Y=\{2\}$, & $\langle Y^\vee \sqcup (E\setminus Y)\rangle=\langle{1, 2^\vee}\rangle$, & $L\cap \langle{1, 2^\vee}\rangle=0 $; \\
   for $Y=\{1, 2\}$, & $\langle Y^\vee \sqcup (E\setminus Y)\rangle=\langle{1^\vee, 2^\vee}\rangle$, & $L\cap \langle{1^\vee,
                                                                                                       2^\vee}\rangle=0$.
  \end{tabular}

\end{Ex}

\begin{Th}\label{sec:theorem1}
  The mapping $\nu_E$ is a bijection between the set of Lagrangian subspaces {\warn $\mathscr{L}_E$}
  and the set $\mathscr{B}_E$ of binary delta-matroids on the set $E$.
\end{Th}

We split the proof of this theorem into several lemmas.

\begin{Def}
  We  say that a Lagrangian subspace $L$ in $V_E$ is {\em graphic} if for each
$e\in E$ there exists an element $v_e\in L$ such that  $(v_e, e)=1$ and $(v_e, e')=0$ for all $e'\in E$, $e'\not= e$.
\end{Def}

By dimension consideration, the collection of such elements $\{v_e\}, e\in E$ forms a basis in the space $L$.

\begin{Ex}% [Mapping $\nu_E$]
  % Let $E$ be a arbitrary Lagrangian subspace in $V_E$.
  % Denote by $\nu_E(L)$ the set system $(E; \Psi_L)$, where the subset $Y\subset E$ belongs to $\Psi_L$ if and only if
  % $L\cap\langle Y^\vee\sqcup(E\setminus )\rangle$

  The Lagrangian subspace $L$ from Example 2.2 is not a graphic one.
  Indeed, for the element $e=1 \in E$, there are two elements $v_e$ such that $(e, v_e)=1$.
  (namely, $1^\vee+2+2^\vee$ and $1+1^\vee+2^\vee$), but for any such element $v_e$ the equality $(2, v_e)=1$ holds as well.

  The subspace $\langle 1^\vee, 2^\vee \rangle$ is an example of a graphic Lagrangian subspace
  in $V_{\langle 1^\vee, 2^\vee \rangle}$. (For $e=1$, we can take $v_e=1^\vee$, for $e=2$ we take $v_e=2^\vee$).
\end{Ex}

\begin{Le}\label{sec:1.01}
  The mapping $\nu_E$ determines a bijection between graphic Lagrangian subspaces in $V_E$ and
  non-degeneracy delta-matroids of framed graphs {\warn on the set of vertex $E$}.
\end{Le}
% ================================================================
\begin{Proof}
  Let $L \subset V_E$ be a graphic Lagrangian subspace;
  assign a symmetric $|E|\times|E^\vee|$-matrix $A(L)$ over $\mathbb{F}_2$
  to this subspace as follows: {put $(v_e, e'^\vee)$
  on the intersection of the row $e$ and
  ({\warn The symmetry of the matrix follows from the fact that~$L$ is Lagrangian:}
    indeed, % for $e\not=e'^\vee$
    the equations $(v_e,e)=(v_{e'},e')=1$ % $(v_e, v_{e^\vee})=(v_e, v_{e'^\vee})=1$
    (for $e\not=e'^\vee$),  $(v_e,e')=(v_{e'},e)=0$ % $(v_e,v_e')=(v_{e'},v_e)=0$
    and  $(v_e,v_{e'})=0$ imply that
    $(v_{e},e'^\vee)=(v_{e'},e^\vee)$ for all $e$ and $e'$).
    One can obtain an arbitrary symmetric matrix in this way.
    Conversely, from a symmetric matrix one can reconstruct the Lagrangian subspace.
    Indeed, $L$ is the Lagrangian subspace in $V_E$ spanned by the vectors $v_e=e^\vee+\sum_{e'\in E}A(L)_{e,e'^\vee}e'$.

    On the other hand, to each framed graph~$G$ with the vertex set~$E$ its adjacency matrix
    $A(G)$ over $\mathbb{F}_2$ is associated.
    By putting $A(L)=A(G)$, we get a one-to-one correspondence between the two sets.
    Let us prove that under this correspondence the set system $\nu_E(L)$ assigned to the Lagrangian subspace~$L$,
    is taken to the non-degeneracy delta-matroid of the graph~$G$.
    Indeed, the subset $Y\subset E$ is feasible, $Y\in \Phi_L$, if and only if the
    sub-matrix $A|_Y$ is non-degenerate over~$\mathbb{F}_2$. The last statement is equivalent to the assertion that the
    subspace $L{\cap} \langle Y^\vee{\sqcup} (E{\setminus} Y)\rangle$ contains only a zero vector.

    Let us prove the last statement.
    The subspace $L{\cap} \langle Y^\vee{\sqcup} (E{\setminus} Y)\rangle$ contains a non-zero vector
    if and only if there exists a non-zero  linear combination
    $\sum_{e\in E}\lambda_ev_e$ (here $v_e=e^\vee+\sum_{e'\in E}A(L)_{e,e'^\vee}e'$) in $L$
    belonging to $\langle Y^\vee{\sqcup} (E{\setminus} Y)\rangle$.
    This means that there exist  $\lambda_e\in\F_2$, $e\in E$, not all equal to~$0$ and such that $\sum_{e\in E}\lambda_ev^*_e=0$,
    where
    \[v^*_e=\left\{
      \begin{aligned}
        \textstyle
        e^\vee+\sum_{e'\in Y}A(L)_{e,e'^\vee}e', &\text{\ if  $e\in E{\setminus}Y$}\\
        \textstyle
        \sum_{e'\in Y}A(L)_{e,e'^\vee}e', &\text{\ if $e\in Y$}\\
      \end{aligned}
    \right.\]
  (here $v^*_e$ is the restriction of $v_e$ to $Y{\sqcup} (E^\vee{\setminus}Y)$).
  This statement is equivalent to degeneracy of the matrix
  \[\begin{pmatrix}
      0 & A|_{Y} \\
      E & *
    \end{pmatrix},\]
}
(here $0$ is the zero matrix of the appropriate size), and hence of the matrix $A|_Y$.  We arrive at a contradiction.
\end{Proof}

For an arbitrary $L\in \cL_E$ and for an arbitrary $e\in E$ denote by $L*e$ the Lagrangian
subspace obtained from $L$ by the linear transformation of {\warn the} space $V_E$ of the form
$e\mapsto e^\vee,\ e^\vee\mapsto e$, acting trivially on the other vectors of the basis.

\begin{Le}\label{sec:1.02}
  For an arbitrary $L\in \cL_E$ and an arbitrary $e\in E$ the following statement is true: \(\nu_E(L)*e=\nu_E(L*e)\).
  In other words, local duality of Lagrangian subspaces descends to twisting of delta-matroids
  under the map $\nu_E$.
\end{Le}

\begin{Proof}
  Let $Y\subset E$ be an arbitrary subset. Note that
  \[
    (L*e){\cap} \langle Y^\vee {\sqcup}(E{\setminus} Y)\rangle=L{\cap} \langle(Y^\vee \Delta \{e^\vee\}){\sqcup}(E{\setminus} (Y \Delta \{e\}))\rangle.
  \]
  It follows that $Y$ is a feasible subset for $\nu_E(L*e)$ if and only if
  $L{\cap} \langle(Y^\vee \Delta \{e^\vee\}){\sqcup}(E{\setminus} (Y \Delta \{e\}))\rangle=0$.
  Thus $Y\Delta e$ is feasible for $\nu_E(L)$ or, equivalently, $Y$ is feasible for $\nu_E(L)*e$.
\end{Proof}

Clearly, the operations  $*e$ and $*e'$ specified by  (not necessarily distinct) elements $e,e'\in E$ commute with each other;
therefore, the operation $*E'$ is well defined for
an arbitrary subset $E'\subset E$.

\begin{Le}\label{sec:1.03}
  For any Lagrangian subspace $L\in \cL_E$,
  there exists a subset $E'\subset E$ such that
  the Lagrangian subspace $L*E'$ is graphic.
\end{Le}
\begin{Proof}
  We start with the choice of a ``good'' basis of $L$. We proceed as follows.

  Choose a vector $e_1$ from the standard basis $E{\sqcup} E^\vee$ of $V_E$ such that there exists a vector $v_1\in L$
  such that $(e_1, v_1)=1$. (Pick $v_1$ for the first element of the ``good basis'').
  Then pick a vector $e_2$ from the standard basis in $V_E$  such that there exists a vector $v_2\in L$,
  with $(e_2, v_2)=1$.
  Add the vector $v'_2=v_2-(e_1,  v_2)v_1$ to the ``good basis''.
  Repeat the procedure
 to obtain a basis in $L$ (similarly to the Gram--Schmidt process).
  Then apply to~$L$ the local duality through the set of those $e_1,e_2,\dots, e_{|E|}$ that belong to $E^\vee$.
  We obtain the subspace $L_1$.
  It  corresponds to the matrix $A(L_1)$ (which is symmetric as long as $L_1$ is a Lagrangian space).
\end{Proof}

\begin{corollary}(follows from Lemmas~\ref{sec:1.01},~\ref{sec:1.02} and ~\ref{sec:1.03} )
  The mapping~$\nu_E$ takes every Lagrangian subspace in $V_E$ to a
  binary delta-matroid over the set~$E$.
\end{corollary}

Now we can complete the proof of Theorem~\ref{sec:theorem1}.

Let us prove that $\nu_E:\cL_E\rightarrow \cB_E$ is an injection.
Suppose the converse. Then there exist distinct Lagrangian subspaces $L_1,L_2\in\cL_E$,
such that $\nu_E(L_1)=\nu_E(L_2)$.
Let $E'\subset E$ be the set corresponding to $L_1$ in Lemma~\ref{sec:1.03}.
Then
\[\nu_E(L_1*E')=\nu_E(L_1)*E'=\nu_E(L_2)*E'=\nu_E(L_2*E'),\]
by Lemma~\ref{sec:1.02}.
But it is shown in Lemma~\ref{sec:1.01}, that the equation  $\nu_E(L_1*E')=\nu_E(L_2*E')$  implies that $L_1*E'=L_2*E'$.
Therefore, $L_1*E'*E'=L_2*E'*E'$, i.e. $L_1=L_2$.

Now let us prove that
$\nu_E:L_E\rightarrow B_E$ is a surjection.
Indeed, for every binary delta-matroid $B\in \cB_E$ there exists a
subset $E'\subset E$ such that $B*E'$ is a graphic delta-matroid.
There exists a Lagrangian subspace $L\in \cL_E$ such that $\nu_E(L)=B*E'$.
Now $\nu_E(L)*E'=B$ and, by Lemma~\ref{sec:1.02}, $\nu_E(L)*E'=\nu_E(L*E')$, i.e. $\nu_E(L*E')=B$.

Theorem~\ref{sec:theorem1} is proven.

\section{Lagrangian subspaces and binary delta-matroids of embedded graphs}
\label{sec:--3.0}

Denote by $\cG_E$ the set of connected ribbon graphs
with the set of ribbons labeled by the elements of~$E$.

In~\cite{BGBS}, a mapping from $\cG_E$ to $\cL_E$ is constructed.
It has the following form.
Let $\Gamma$ be a connected ribbon graph with the set of ribbons~$E$
interpreted as the union of two sets of closed topological disks called {\em vertices}
$V(G)$ and {\em edges} $E(G)$ satisfying the following conditions:
% \begin{minipage}{1.0\linewidth}
  \begin{itemize}
    \setlength\itemsep{-5pt}
  \item edges and vertices intersect by disjoint line segments;
  \item each such segment lies in the closure of precisely one edge and one vertex;
  \item each edge contains two such segments.
  \end{itemize}
% \end{minipage}

%% ================================================================
Given a ribbon graph~$\Gamma$, remove small open discs from the centers of
the vertices, which are discs. Let ~$F_\Gamma$ denote the resulting  two-dimensional surface with a boundary.

To each $e\in E$, we associate~$h_e$, an element of the relative homology group  $H_1(F_\Gamma,\partial F_\Gamma)$.
This element is represented by a segment  going along the edge~$e$ and connecting the boundaries of the discs
that are removed from the vertices incident to the edge $e$).

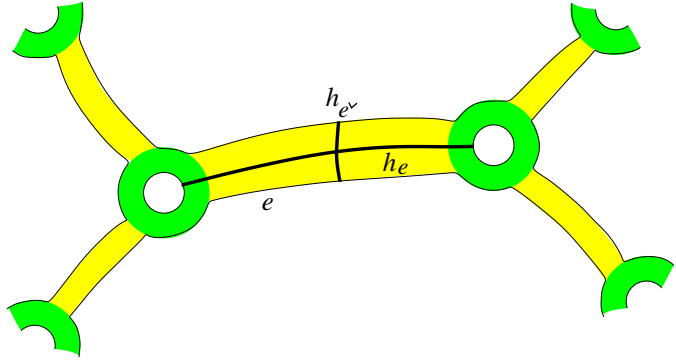
\begin{figure}
\centering
\resizebox{.75\textwidth}{!}{\input{ris.pspdftex}}
  \caption{A ribbon graph without discs removed around the centers of the vertices, with elements ~$h_e$,~$h_{e^\vee}$
    of the first relative homology group $H_1(F_\Gamma,\partial F_\Gamma)$ assigned to the edge~$e$}\label{fff1}
\end{figure}

On the other side, to each element  $e^\vee\in E^\vee$ we may associate an element  $h_{e^\vee}$
in the relative homology group $H_1(F_\Gamma,\partial F_\Gamma)$
that is represented by a segment that goes across the edge~$e$ and connects the opposite sides of this
edge  (see Fig.~\ref{fff1}).

To each continuous cycle $\gamma: S^1\to F_\Gamma$, we associate the vector
$\sum_{e\in E}((\gamma,h_e)h_e+(\gamma,h_{e^\vee})h_{e^\vee})$ in ~$V_E$.
(The brackets $(\cdot,\cdot)$ in this formula denote the intersection form between
the first absolute and relative homology for the given surface with boundary $F_\Gamma$).
%Notice that the space $V_E$ can be naturally presented as the vector
%subspace in the first relative homology  $H_1(F_\Gamma,\partial F_\Gamma)$
%with coefficients in ~$\F_2$, with a basis of~$V_E$ formed by the vectors $h_e,h_{e^\vee}$.
%Notice that the intersection form defines a symplectic structure in $V_E$. This symplectic structure coincides
%with the one introduced in Sec.~\ref{sec:binary-delta-matr}.
As shown in~\cite{BGBS,KS}, the subspace of~$V_E$ formed by the vectors that correspond
to all cycles~$\gamma$, is Lagrangian.  Denote this subspace by~$\pi_E(\Gamma)$.

On the other hand, Bouchet~\cite{Bo89} assigned to each ribbon graph
a set system whose ground set is the set of edges of the graph:
a subset of edges is feasible if the restriction of the given graph to this subset is a {\it quasi-tree},
that is, a ribbon graph with a connected boundary. Bouchet showed that the set system assigned
to a ribbon graph in such a way is a delta-matroid. We denote this delta-matroid by~$\rho_E(\Gamma)$.

\begin{Th}
  The mapping~$\nu_E$ is compatible with the mappings~$\pi_E$ and $\rho_E$.
  Namely, for an arbitrary  $\Gamma\in \cG_E$ the following identity holds:\\ $\rho_E(\Gamma)=\nu_E(\pi_E(\Gamma))$.
\end{Th}
\begin{Proof}
  Let first $\Gamma$ be a ribbon graph with a single vertex, i.e. a (framed) chord diagram.
   Then the statement is true, since both mappings are compatible with the mapping which assigns to a chord diagram $\Gamma$
    the adjacency matrix of its intersection graph.
  Conversely, each of the mappings is compatible with the twist operation on the corresponding ribbon graphs
  $\rho_E(\Gamma*e)=(\pi_E(\Gamma))*e$.
  For an arbitrary ribbon graph $\Gamma$ find a set  $E'\subset E$ such that  $\Gamma*E'$
  has a single vertex; then
  $(\rho_E(\Gamma))*E'=\nu_E(\pi_E(\Gamma)*E')=\nu_E(\pi_E(\Gamma))*E'$,
  i.e. $(\rho_E(\Gamma))*E'=\nu_E(\pi_E(\Gamma))*E'$, hence
  $\rho_E(\Gamma)=\nu_E(\pi_E(\Gamma))$ as required.
\end{Proof}

\section{Hopf Algebras Isomorphism}
\label{sec:1}

Let  $n=|E|$.
Denote by $\cL_n$ the set of isomorphism classes of Lagrangian subspaces
$\cL_E\subset V_E$ with respect to bijections of $n$-element sets.

Let $\cB_n$ denote the set of isomorphism classes of binary delta-matroids on~$n$ elements.

Klepsyn and Smirnov in \cite{KS} introduce the structure of a graded commutative and cocommutative Hopf algebra
on the infinitely dimensional vector space
\[
\mathbb{C}\cL =
\mathbb{C}\cL_0\oplus\mathbb{C}\cL_1\oplus \cdots,
\]
where $\mathbb{C}\cL_n$ is the vector space over $\mathbb{C}$ freely spanned by the set $\cL_n$.
 Multiplication in this Hopf algebra is given by the operation of direct sum of
 Lagrangian subspaces in the direct sum of symplectic spaces, which is extended to~$\mathbb{C}\cL$  by linearity.
The comultiplication $\mathbb{C}\cL\to\mathbb{C}\cL\otimes\mathbb{C}\cL$
assigns to a Lagrangian subspace $L\subset V_E$
the sum of the tensor products of the Lagrangian subspaces
\[
L\mapsto \sum_{I\subset E} L_I\otimes L_{E{\setminus} I},
\]
where, for a subset~$I$ of the set ~$E$, ~$L_I\subset V_I$  denotes the subspace,
which is the symplectic reduction of the Lagrangian subspace~$L$ (see ~\cite{KS}).
This multiplication can be naturally transferred to the vector space~$\mathbb{C}\cL$,
spanned by the Lagrangian subspaces, considered up to renumbering finite element sets.

Meanwhile, in~\cite{LZ}, a graded Hopf algebra of binary delta-matroids is constructed
\[
\mathbb{C}\cB = \mathbb{C}\cB_0\oplus\mathbb{C}\cB_1\oplus \cdots,
\]
where the subspace $\mathbb{C}\cB_n$
is freely spanned over $\mathbb{C}$ by the set $\cB_n$.
The multiplication in this Hopf algebra is given by the direct sum of set systems
extended to~$\mathbb{C}\cB$ by linearity. The coproduct of a given set system $(E;\Psi)$
is the sum
\[
\mu(E;\Psi)=\sum_{E'\subset E}\Psi|_{E'}\otimes\Psi|_{E{\setminus} E'},
\]
where the set $\Psi|_{E'}$ consists of those elements of the set ~$\Psi$
that are contained in ~$E'$.

The mapping ~$\nu_E$ (see Def.~\ref{sec:def1}) is equivariant
with respect to bijections of finite sets both on the set of Langrangian subspaces and on the set of binary delta-matroids.
Hence the set of such mappings defines a graded linear mapping

\[
\nu:\mathbb{C}\cL\to \mathbb{C}\cB,\qquad \nu_n:\mathbb{C}\cL_n\to \mathbb{C}\cB_n, \qquad n=0,1,2,\dots.
\]

This linear mapping appears to be an isomorphism:

\begin{Th}
 The mapping~$\nu:\mathbb{C}\cL\to \mathbb{C}\cB$ is a graded isomorphism of Hopf algebras.
\end{Th}
\begin{Proof}
 The mapping~$\nu$ transfers the multiplication and the comultiplication
 in the Hopf algebra of Lagrangian subspaces to the multiplication and the comultiplication,
 respectively, in the algebra of binary delta-matroids. This can be seen from the definitions above.
\end{Proof}

\section{Four-term relations and weight systems}

In~\cite{V} V.~A.~Vassiliev introduced the our-term relations for functions on chord
diagrams. He proved that any invariant of order at most~$n$
determines a function on chord diagrams that satisfies these relations.
Such a function is called a  {\it weight system}. Every four-term relation
corresponds to a chord diagram and to a pair of chords with neighboring ends in it.
The remaining three diagrams that participate in this relation
can be built from the initial one by application of one of the two
(mutually commuting) Vassiliev moves, and their compositions.
 In ~\cite{L} Vassiliev moves were extended to framed diagrams,
 which are chord diagrams associated to ribbon graphs with possibly twisted ribbons, and
 the corresponding four-term relations were described.

Kleptsyn and Smirnov in \cite{KS} extended Vassiliev moves to Lagrangian subspaces.
Let, as above, $E$ be a finite set, $V_E$ be the vector space over~$\F_2$
spanned by the elements of the set $E{\sqcup} E^\vee$, and let $e,e'\in E$ be two distinct elements in ~$E$.
Then the {\it first Vassiliev move\/}, assigned to a pair~$e,e'$, is a linear mapping
$V_E\to V_E$ preserving all the basis vectors except for the vectors $e^\vee,e'^\vee$.
The action on these vectors is defined as follows:
\[
e^\vee\mapsto e^\vee+e';\qquad e'^\vee\mapsto e'^\vee+e.
\]
Notice that  the first Vassiliev move is symmetric with respect to  the transposition of the elements  ~$e$ and ~$e'$.

 The {\it second Vassiliev move\/} for the pair ~$e,e'$ is a linear mapping
$V_E\to V_E$ obtained from the first move by conjugation with respect to the twist along the element $e'\in E$,
see Sec.~\ref{sec:binary-delta-matr}. In contrast to the first move,
the description of the second one depends on the order of elements in the pair ~$e,e'$.
The action of each Vassiliev move
on the set of Lagrangian subspaces is induced by its action on ~$V_E$.

In  ~\cite{LZ}, the authors define the first and the second Vassiliev moves for binary delta-matroids  $\cB_E$.
To define the second Vassiliev move,
they use the recently introduced (see ~\cite{MMB}) concept of  handle sliding for delta-matroids.
In ~\cite{LZ}, it is shown (see Proposition  4.10) that the action of the first and second Vassiliev moves on the space ~$V_E$
 as defined by Kleptsyn--Smirnov  coincides with the one defined by Zhukov and Lando for binary delta-matroids.
 Taking into account Theorem ~\ref{sec:theorem1}, we obtain the following statement.

\begin{Th}
The graded Hopf algebras isomorphism~$\nu:\mathbb{C}\cL\to \mathbb{C}\cB$ descends to a graded quotient Hopf algebras isomorphism $\nu:\cF\mathbb{C}\cL\to \cF\mathbb{C}\cB$, that of  the  Hopf algebras $\mathbb{C}\cL$ and  $\mathbb{C}\cB$ modulo the corresponding four-term relations.
\end{Th}

\end{document}

%% file: ris.pspdftex
\begin{picture}(0,0)%
\includegraphics{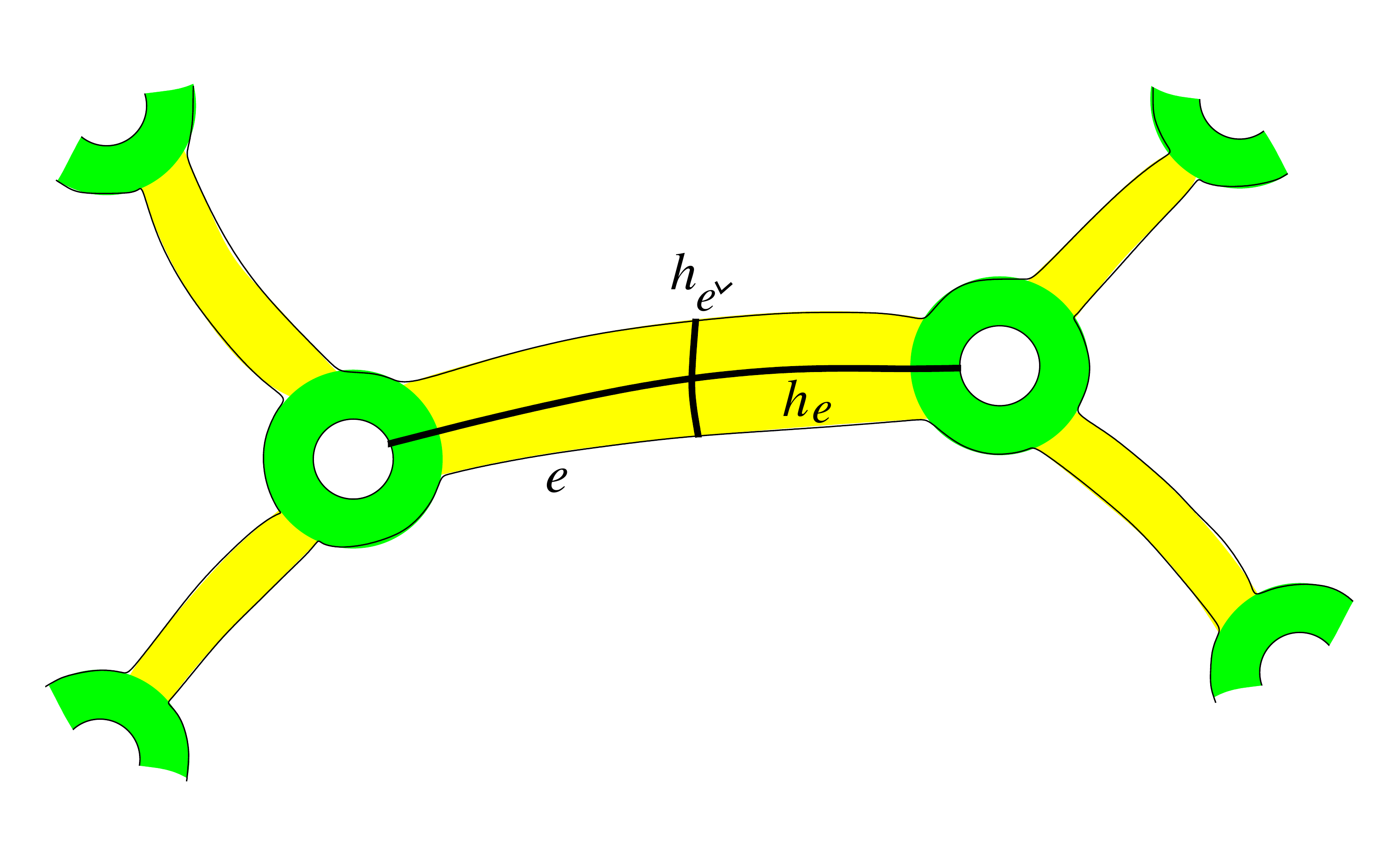}%
\end{picture}%
\setlength{\unitlength}{3947sp}%
\begingroup\makeatletter\ifx\SetFigFont\undefined%
\gdef\SetFigFont#1#2#3#4#5{%
  \reset@font\fontsize{#1}{#2pt}%
  \fontfamily{#3}\fontseries{#4}\fontshape{#5}%
  \selectfont}%
\fi\endgroup%
\begin{picture}(15750,9645)(76,-9871)
\end{picture}%

%% file: Zhukov_eng18Arxiv.bbl
\begin{thebibliography}{99}
\bibitem{BGBS} Booth, Richard F.; Borovik, Alexandre V.; Gelfand, Israel M.; Stone,
  David A. {\it Lagrangian matroids and cohomology}. Ann. Comb. 4 (2000), no. 2,
  171–182.
  % \bibitem{BGBS} Richard F. Booth, Israel M. Gelfand,
  %   Alexandre V. Borovik, David A. Stone,
  %   {\it Lagrangian Matroids and Cohomology},
  %   Journal of Algebraic Combinatorics 8 (1998), 235–252

\bibitem{Bo871} A.Bouchet,
  {\it Greedy algorithms and symmetric matroids},
  Math. Programm. 38 (1987), 147--159

\bibitem{Bo872} Bouchet, A.(F-LMNS-IC)
  {\it Representability of $\Delta$-matroids}. Combinatorics (Eger, 1987), 167–182,
  Colloq. Math. Soc. Janos Bolyai, 52, North-Holland, Amsterdam, 1988.
  % \bibitem{Bo872} A.Bouchet,
  %   {\it Representability of $\Delta$-matroids},
  %   in: Proceedings of the 6th Hungarian Colloquium of Combinatorics,
  %   Colloq. Math. Soc. J\'anos Bolyai 38 (1987), 167--182

\bibitem{Bo89} A.Bouchet,
  {\it Maps and delta-matroids},
  Discrete Math. 78 (1989), 59--71


\bibitem{C09} S.~Chmutov,
  {\it Generalized duality for graphs on surfaces and the signed Bollob\'as--Riordan
    polynomial},
  J. of Combin. Theory Ser. B 99 (2009) 617--638


\bibitem{CMNR} C.~Chun, I.~Moffatt, S.~D.~Noble, R.~Rueckriemen,
  {\it Matroids, delta-matroids and embedded graphs},
  arXiv: 1403.0920v1, 45 pp.

\bibitem{KS} V.~Kleptsyn, E.~Smirnov
  {\it Ribbon graphs and bialgebra of Lagrangian subspaces},
  Journal of Knot Theory and Its Ramifications Vol. 25, No. .12 (2016) 1642006

  % \bibitem{KS} V.~Kleptsyn, E.~Smirnov
  %   {\it Ribbon graphs and bialgebra of Lagrangian subspaces}, arxiv:1401.6160

\bibitem{L00}      {S.~K.~Lando,}
  {\it On a Hopf algebra in graph theory},
  J. Comb. Theory, Ser. B, vol.~{\bf 80} (2000), 104--121.

\bibitem{L}      {S.~K.~Lando,}
  {\it $J$-invariants of ornaments and framed chord diagrams},
  Funct. Anal. Appl.,~{\bf 40}(1) (2006), 1--13.


\bibitem{LZ} S.~Lando, V.~Zhukov,
  {\it Delta-matroids and Vassiliev invariants}, arxiv:1602.00027

\bibitem{LZ04} S.~Lando, A.~Zvonkin, {\it Graphs on surfaces and their applications}, Springer, 2004.

\bibitem{MG} Malic G.
  {\it An action of the Coxeter group BCn on maps on surfaces,
    Lagrangian matroids and their representations}
  arXiv:1507.01957v3
  %
  % \bibitem{MG} Malic G. {\it Representations of
  %   Lagrangian matroids associated to partial duals
  %   of maps on surfaces} arXiv preprint
  %   arXiv:1507.01957 (2015).

\bibitem{MMB} Iain Moffatt, Eunice Mphako-Banda,
  {\it Handle slides for delta-matroids}, arXiv:1510.07224, 12 pp.


\bibitem{V} V.~A.~Vassiliev, {\it Cohomology of knot spaces}, in: Theory of singularities
and its applications, 23-69, Adv. Soviet Math., 1, Amer. Math. Soc.,
Providence, RI, 1990.

%% {\TODO
%% \bibitem{VW} Andrew Vince, Neil White {\it Orthogonal Matroids}
%% }

\end{thebibliography}
